\documentclass[12pt]{article}
\usepackage{latexsym, amssymb}
\usepackage{mathbbol}

\makeatletter
\def\eqnarray{\stepcounter{equation}\let\@currentlabel=\theequation
\global\@eqnswtrue
\tabskip\@centering\let\\=\@eqncr
$$\halign to \displaywidth\bgroup\hfil\global\@eqcnt\z@
  $\displaystyle\tabskip\z@{##}$&\global\@eqcnt\@ne
  \hfil$\displaystyle{{}##{}}$\hfil
  &\global\@eqcnt\tw@ $\displaystyle{##}$\hfil
  \tabskip\@centering&\llap{##}\tabskip\z@\cr}

\def\endeqnarray{\@@eqncr\egroup
      \global\advance\c@equation\m@ne$$\global\@ignoretrue}
\makeatother

\begin{document}
\bibliographystyle{tom}

\newtheorem{lemma}{Lemma}[section]
\newtheorem{thm}[lemma]{Theorem}
\newtheorem{cor}[lemma]{Corollary}
\newtheorem{voorb}[lemma]{Example}
\newtheorem{rem}[lemma]{Remark}
\newtheorem{prop}[lemma]{Proposition}
\newtheorem{stat}[lemma]{{\hspace{-5pt}}}

\newenvironment{remarkn}{\begin{rem} \rm}{\end{rem}}
\newenvironment{exam}{\begin{voorb} \rm}{\end{voorb}}

\newcounter{teller}
\renewcommand{\theteller}{\Roman{teller}}
\newenvironment{tabel}{\begin{list}%
{\rm \bf \Roman{teller}.\hfill}{\usecounter{teller} \leftmargin=1.1cm
\labelwidth=1.1cm \labelsep=0cm \parsep=0cm}
                      }{\end{list}}

\newcommand{\Ni}{{\bf N}}
\newcommand{\Ri}{{\bf R}}
\newcommand{\Ci}{{\bf C}}

\newcommand{\proof}{\mbox{\bf Proof} \hspace{5pt}} 
\newcommand{\remark}{\mbox{\bf Remark} \hspace{5pt}}
\newcommand{\ruimte}{\vskip10.0pt plus 4.0pt minus 6.0pt}

\newcommand{\RRe}{\mathop{\rm Re}}
\newcommand{\supp}{\mathop{\rm supp}}
\newcommand{\one}{\mathbb{1}}

\hyphenation{groups}
\hyphenation{unitary}

\thispagestyle{empty}

\begin{center}
{\Large\bf Contraction semigroups on $L_\infty(\Ri)$ } \\[5mm]
{\large A.F.M. ter Elst$^1$ and Derek W. Robinson$^2$ } \\[10mm]
\bf Dedicated to the memory of G\"unter Lumer 1929--2005
\end{center}

\vspace{5mm}

\begin{center}
{\bf Abstract}
\end{center}

\begin{list}{}{\leftmargin=1.8cm \rightmargin=1.8cm \listparindent=10mm 
   \parsep=0pt}
\item

If $X$ is a non-degenerate vector field on $\Ri$  and $H=-X^2$ we examine conditions 
for the closure of $H$ to generate a continuous semigroup on $L_\infty$ which extends
to the $L_p$-spaces.
We give an example which cannot be extended and an example which extends but 
for which the real part of the generator on $L_2$ is not lower semibounded.
\end{list}

\vspace{7cm}
\noindent
June 2006

\vspace{5mm}
\noindent
AMS Subject Classification: 47B44, 58G03.

\vspace{5mm}

\noindent
{\bf Home institutions:}    \\[3mm]
\begin{tabular}{@{}cl@{\hspace{10mm}}cl}
1. & Department of Mathematics  & 
  2. & Centre for Mathematics   \\
& University of Auckland   & 
  & \hspace{15mm} and its Applications  \\
& Private bag 92019 & 
  & Mathematical Sciences Institute  \\
& Auckland & 
  & Australian National University  \\
& New Zealand  & 
  & Canberra, ACT 0200  \\
& & 
  & Australia  \\[8mm]
\end{tabular}

\newpage
\setcounter{page}{1}

\section{Introduction}\label{Svf1}

The Lumer--Phillips theorem \cite{LuP} is a cornerstone of the theory of continuous semigroups.
The theorem characterizes the generator of a contraction semigroup with the aid of a dissipativity
condition. 
The latter is based on the elementary properties of the operator $-d^2/dx^2$ of double differentiation
acting on $C_0(\Ri)$.
In this note we analyze  contraction semigroups $S$ generated by squares $-X^2$ of vector fields
$X=a\,d/dx$ acting on $C_0(\Ri)$, or $L_\infty(\Ri)$.
An integral part of the analysis consists of examining the one-parameter groups $T$ generated by $X$.
Throughout we assume $a>0$.
If $a$ is smooth this is the one-dimensional analogue of  
H\"ormander's condition \cite{Hor1}.

First, we identify the kernel of $S$ acting on $L_\infty(\Ri)$.
Secondly, $T$ is defined as a weak$^*$ continuous group of contractions on $L_\infty$ and 
we derive necessary and sufficient conditions for it to extend to a continuous group on the
$L_p(\Ri\,;\rho\,dx)$-spaces with $p \in [1,\infty\rangle$, where 
$\rho\colon \Ri\to \langle0,\infty\rangle$ is a $C^\infty$-function.
These conditions also ensure that $S$ extends to a continuous semigroup.
Thirdly, we characterize those $S$, or $T$, which extend to a contraction semigroup, or group,  on 
$L_p(\Ri\,;\rho\,dx)$ for some $p \in [1,\infty\rangle$.
Fourthly, we give an example of a smooth vector field with a uniformly bounded coefficient
for which neither $T$ nor $S$ can be extended to any of the $L_p$-spaces with $p<\infty$.
Fifthly, we give an example of a smooth vector field with a uniformly bounded coefficient
which is uniformly bounded away from zero for which $T$ and $S$ extend to all the $L_p$-spaces but 
the real part of the generator of $S$ on $L_2(\Ri\,;\rho\,dx)$ is not lower semibounded.
In particular the $L_2$-generator cannot satisfy a G\aa rding inequality.
Since the G\aa rding inequality is the usual starting point for the analysis of elliptic
divergence form operators  on $L_2(\Ri\,;\rho\,dx)$, e.g., operators of the form $X^*X$, this example
clearly demonstrates that the theory of `non-divergent' form operators such as $-X^2$ on $L_\infty(\Ri)$
is very different.
Finally we discuss the volume doubling property for balls (intervals) whose radius (length)
is measured by the distance associated with $X$.

\section{Preliminaries}\label{Svf2}

Let $a\colon \Ri\to \langle0,\infty\rangle$ be a locally bounded differentiable function
and assume the derivative $a'$ is locally bounded.
Further assume 
\begin{equation}
\int^\infty_0dx\,a(x)^{-1}=\infty=\int^0_{-\infty}dx\,a(x)^{-1}
\label{edvf1.1}
\;\;\;.
\end{equation}
Equip $\Ri$ with the measure $\rho\, dx$ where  
$\rho\colon \Ri\to \langle0,\infty\rangle$ is a $C^\infty$-function.
Consider the vector field $X=a\,d/dx$ and the corresponding operators $X_{\rm min}$ and 
$X_{\rm max}$ on $L_\infty(\Ri\,;\rho \,dx)$ with domains
$D(X_{\rm min})=C_c^\infty(\Ri)$ and $D(X_{\rm max}) = C_c^1(\Ri)$.
Set $H_{\rm min} = - X_{\rm min}^{\,2}$ and $H_{\rm max} = - X_{\rm max}^{\,2}$.
Since we are dealing with operators on $L_\infty$ it is appropriate to deal with the weak$^*$ topology.

\begin{prop}\label{pdvf1.1}
\mbox{}
\begin{tabel}
\item\label{pdvf1.1-1}
The operators $X_{\rm min}$ and 
$X_{\rm max}$ are weak$^*$ closable
and $\overline X_{\rm min}=\overline  X_{\rm max}$,
where the bar denotes the weak$^*$ closure.
\item\label{pdvf1.1-2}
The operator $H_{\rm max}$ is weak$^*$ closable and its weak$^*$ closure 
$\overline H_{\rm max}$ generates a 
semigroup $S$ which 
is weak$^*$ continuous, positive, contractive and holomorphic in the 
open right half-plane.
\item\label{pdvf1.1-3}
$\overline H_{\rm max}= -  \overline X_{\rm max}^{\;2} $
and in particular $\overline X_{\rm max}^{\;2} $ is weak$^*$ closed.
\item\label{pdvf1.1-4}
If $a \in C^\infty(\Ri)$ then  $\overline H_{\rm min} = \overline H_{\rm max}$, where
 $\overline H_{\rm min}$ is the weak$^*$ closure of $H_{\rm min}$.
\end{tabel}
\end{prop}
\proof\
For all $x_0\in \Ri$ the ordinary differential equation $\dot{x}=a(x)$,  with initial data
$x(0)=x_0$, has a unique maximal solution which we denote by $t\mapsto e^{tX}x_0$.
Since $a$ satisfies~(\ref{edvf1.1}) this maximal solution is defined for all $t\in\Ri$.
Moreover, $e^{sX}e^{tX}x_0=e^{(s+t)X}x_0$ and 
\begin{equation}
\int_{x_0}^{e^{tX}x_0}dx\,a(x)^{-1}=t
\label{edvf1.2}
\end{equation}
for all $s,t\in\Ri$ and $x_0\in \Ri$.
In addition both the maps $t\mapsto e^{tX}x_0$ and $x\mapsto e^{sX}x$ are continuous.
In particular for all $t \in \Ri$ the map  $T_t \colon L_\infty \to L_\infty$ 
defined by $(T_t\varphi)(y)=\varphi(e^{-tX}y)$ is an isometry and 
$T$ is a  
weak$^*$ continuous group on $L_\infty$.
This group is automatically positive and we next show that 
its generator is the weak$^*$ closure 
of the operator $X_{\rm min}$ on $L_\infty$.

Clearly $X_{\rm min} \subseteq X_{\rm max}$ and by a standard regularization argument 
it follows that $\overline X_{\rm min} = \overline X_{\rm max}$.
Hence to simplify notation we now set 
$X_0=\overline X_{\rm min}= \overline X_{\rm max}$.

One computes from (\ref{edvf1.2}) that
\[
{{d}\over{dy}}\,e^{tX}y={{a(e^{tX}y)}\over{a(y)}}
\]
for all $t\in \Ri$ and $y\in \Ri$. 
Therefore 
\[
\frac{d}{dy} (T_t \varphi)(y) = \varphi'(e^{-tX}y) \cdot \frac{a(e^{tX}y)}{a(y)}
\]
for all $\varphi \in D(X_{\rm max})$, $y \in \Ri$ and $t > 0$.
So $T_t (D(X_{\rm max})) \subseteq D(X_{\rm max})$ for all $t > 0$.
Moreover, 
\begin{eqnarray*}
t^{-1} (\varphi - T_t \varphi)(y)
& = & - t^{-1} \int_0^t ds \, \frac{d}{ds} \varphi(e^{-sX}y)  \\[5pt]
& = & t^{-1} \int_0^t ds \, \varphi'(e^{-sX}y) \, a(e^{-sX}y) 
= t^{-1} \int_0^t ds \, (T_sX_{\rm max} \varphi)(y)
\end{eqnarray*}
for all $\varphi \in D(X_{\rm max})$, $t > 0$ and $y \in \Ri$, since $\varphi'$ is continuous.
So $\lim_{t \to 0} t^{-1}(I - T_t) \varphi = X_{\rm max} \varphi$ strongly in $L_\infty$
and $X_{\rm max}$ is the restriction of the generator of $T$.
Since $D(X_{\rm max})$ is invariant under $T$ and 
weak$^*$ dense it follows from Corollary~3.1.7 of \cite{BR1} 
that $X_0=\overline X_{\rm max}$ is the generator of $T$.

Next define the semigroup $S$ by the integral algorithm
\begin{equation}
S_t = (4\pi t)^{-1/2} \int^\infty_{-\infty} ds \, e^{-s^2 (4t)^{-1}} T_s
\;\;\;.
\label{edvf1.3}
\end{equation}
Obviously $S$ is weak$^*$ continuous, positive, contractive
and holomorphic in the open right half-plane.
Let $H_0$ denote the weak$^*$ closed  generator of $S$.
If $\varphi\in D( X_0^{\;2})$ then 
\begin{eqnarray*}
t^{-1}\,(I-S_t)\varphi&=&
t^{-1}\,(4\pi t)^{-1/2} \int^\infty_{-\infty} ds \, e^{-s^2 (4t)^{-1}}(I- T_s)\varphi\\[5pt]
&=&t^{-1}\,(4\pi t)^{-1/2} \int^\infty_{-\infty} ds\, e^{-s^2 (4t)^{-1}}
\int^s_0du \, (s-u)\,T_u\,X_0^{\;2}\varphi\\[5pt]
&=&(4\pi)^{-1/2} \int^\infty_{-\infty} ds\, e^{-s^2 /4}
\int^s_0du \, (s-u)\,T_{t^{1/2}u}\, X_0^{\;2}\varphi
\end{eqnarray*}
and it follows in the weak$^*$ limit $t\to0$ that $\varphi\in D(H_0)$.
Hence
$H_0\supseteq  - X_0^{\;2}$.
To prove $H_0= - X_0^{\;2}$ it suffices to establish  that the range $R(I- X_0^{\;2})$
of $I- X_0^{\;2}$ is equal to $L_\infty$.
But  $X_0$
generates the continuous group $T$.
Therefore  $R(I\pm  X_0)=L_\infty$.
Moreover, $I- X_0^{\;2}=(I- X_0)(I+ X_0)$.
Hence  $R(I- X_0^{\;2})=L_\infty$ and $H_0=  -X_0^{\;2}$.

Clearly $H_{\rm max} \subseteq  - X_0^{\;2}=H_0$ so $H_{\rm max}$ is weak$^*$ closable. 
It remains to prove that the weak$^*$ closure $\overline H_{\rm max}$ of $H_{\rm max}$ 
is equal to $H_0$.

Since $T_t D(X_{\rm max}) \subseteq D(X_{\rm max})$ and $X_{\rm max} T_t \varphi=T_t X_{\rm max}\varphi$ for all 
$\varphi\in D(X_{\rm max})$ one deduces  by iteration that 
$T_t D(X_{\rm max}^{\;2})\subseteq D(X_{\rm max}^{\;2})$ and $X_{\rm max}^{\;2} T_t \varphi = T_t X_{\rm max}^{\;2}\varphi$ for all 
$\varphi\in D(X_{\rm max}^{\;2})$.
Next it follows from (\ref{edvf1.3}), by a Riemann approximation argument, that 
$S_t D(X_{\rm max}^{\;2}) \subseteq D(\overline{X_{\rm max}^{\;2}})$  and
$\overline{X_{\rm max}^{\;2}} S_t\varphi = S_t {X_{\rm max}^{\;2}} \varphi$ for all 
$\varphi\in D({X_{\rm max}^{\;2}})$ and all $t > 0$.
Since $S_t$ is continuous it further follows  that 
$S_t D(\overline{X_{\rm max}^{\;2}}) \subseteq D(\overline{X_{\rm max}^{\;2}})$
for all $t > 0$.
But $C_c^1(\Ri) \subseteq D(X_{\rm max}^{\;2}) \subseteq D(\overline H_{\rm max})$ is weak$^*$ dense
in $L_\infty$ by the assumed differentiability of $a$.
Hence by Corollary~3.1.7 of \cite{BR1} it follows that $D(\overline H_{\rm max})$ 
is a core of $H_0$.
Therefore $\overline H_{\rm max} = H_0$.

Finally, if $a \in C^\infty(\Ri)$ then $C_c^\infty(\Ri)$ is a core for $X_{\rm max}^{\;2}$.
Therefore $\overline H_{\rm min} \supseteq H_{\rm max}$. 
Since $H_{\rm min} \subseteq H_{\rm max}$ this completes the proof of the 
proposition.\hfill$\Box$

\begin{remarkn}\label{rvf1.1}
It follows by definition that $T_t C_0(\Ri)\subseteq C_0(\Ri)$ for all $t\in\Ri$
and a simple estimate shows that the restriction of $T$ to $C_0(\Ri)$ is strongly continuous.
Therefore $S_t C_0(\Ri)\subseteq C_0(\Ri)$ for all $t>0$ and the restriction of $S$ to
$C_0(\Ri)$ is also strongly continuous. 
This is a direct consequence of the algorithm (\ref{edvf1.3}).
Thus $T$ is a Feller group and $S$ is a Feller semigroup.
Now let $X_{00}$ and $H_{00}$ denote the generators of the restricted group and the restricted
semigroup, respectively.
Then a slight modification of the foregoing argument allows  one to obtain similar 
characterizations of the generators but in terms of norm closures.
For example, $X_{00}$ is the norm closure of $X_{\rm min}$ which is equal to the norm closure
of $X_{\rm max}$.
The discussion of $H_{00}$ can in fact be simplified.
Since $X_{00}$ generates a strongly continuous group of isometries the operator $-X_{00}^{\;2}$ is
dissipative in the sense of Lumer and Phillips \cite{LuP} and it is norm closed by standard estimates
(see, for example, \cite{Robm}  Lemma~III.3.3).
But one again has  $R(I\pm X_{00})=L_\infty$.
Therefore  $R(I- X_{00}^{\;2})=L_\infty$.
Then $-X_{00}^{\;2}$ generates a strongly continuous
contraction semigroup by the Lumer--Phillips theorem
and it follows by uniqueness that $H_{00}=-X_{00}^{\;2}$.
\end{remarkn}

One can associate  a distance with the vector field $X$ by the definition
\begin{equation}
d(x\,;y)=\sup\{|\psi(x)-\psi(y)|\,;\,\psi\in C_c^\infty(\Ri)\,,\,\|X\psi\|_\infty\leq 1\,\}
\;\;\;.
\label{edvf1.5}
\end{equation}
Clearly one has
\[
|\psi(x)-\psi(y)|=\Big|\int^y_xdz\,\psi'(z)\Big|\leq \Big|\int^y_xdz\,a(z)^{-1}\Big|
\]
for all $\psi\in C_c^\infty(\Ri)$ with $\|X_{\rm min}\psi\|_\infty\leq 1$.
So
\[
d(x\,;y)\leq \Big|\int^y_xdz\,a(z)^{-1}\Big|
\;\;\;.
\]
But by regularizing $a^{-1}$ on a compact interval one deduces that the inequality is in
fact an equality, i.e.,
\[
d(x\,;y)=\Big|\int^y_xdz\,a(z)^{-1}\Big|
\]
for all $x,y\in \Ri$.
Note that by setting $x=e^{-sX}y$  and using (\ref{edvf1.2}) one finds
\begin{equation}
d(e^{-sX}y\,;y)=\Big|\int^{e^{-sX}y}_ydz\,a(z)^{-1}\Big|=|s|
\;\;\;.\label{edvf1.61}
\end{equation}
Therefore the distance is invariant under the flow in the sense that 
\[
d(e^{-tX}x\,;e^{-tX}y)=d(x\,;y)
\]
for all $x,y\in \Ri$ and all $t\geq 0$.
This follows by setting $x=e^{-sX}y$ and
\[
d(e^{-tX}x\,;e^{-tX}y)= d(e^{-sX}e^{-tX}y\,;e^{-tX}y)=|s|
=d(e^{-sX}y\,;y)=d(x\,;y)
\;\;\; ,  \]
where we have  used (\ref{edvf1.61}).

Now one can  calculate the kernel of the semigroup $S$.

\begin{prop}\label{pdvf1.21}
The kernel $K$ of the semigroup $S$ on $L_\infty(\Ri)$ is given by
\begin{equation}
K_t(x\,;y)
= (4\pi t)^{-1/2} \, (a(y) \rho(y))^{-1} e^{-d(x;y)^2 (4t)^{-1}}
\label{edvf1.9}
\end{equation}
for all $x,y\in \Ri$ and $t>0$.
Moreover, $K_t$ is continuous and $\int dy \, \rho(y) \, K_t(x\,;y) = 1$ for all $x \in \Ri$.
\end{prop}
\proof\
First by (\ref{edvf1.3}) one has 
\[
(S_t\varphi)(x)=(4\pi t)^{-1/2}\int^\infty_{-\infty}ds\, e^{-s^2 (4t)^{-1}}\varphi(e^{-sX}x)
\]
for all $\varphi\in C_c^\infty(\Ri)$,  $t>0$ and $x\in \Ri$.
Therefore by a change of variables $y=e^{-sX}x$  one deduces that
\[
(S_t\varphi)(x)=(4\pi t)^{-1/2}\int^\infty_{-\infty}dy\,a(y)^{-1}
 e^{-d(x;y)^2 (4t)^{-1}}\varphi(y)
\]
since $|s|=d(x\,;y)$ by (\ref{edvf1.61}).
The representation (\ref{edvf1.9}) follows immediately.

Clearly $K_t$ is continuous and $H_{\rm max} \one = 0$.
So $S_t \one = \one$ in $L_\infty$-sense. 
Therefore $\int dy \, \rho(y) \, K_t(x\,;y) = 1$ for all $t > 0$ and almost every $x \in \Ri$.
Moreover, the map $x \mapsto \int dy \, \rho(y) \, K_t(x\,;y)$ is continuous. 
Hence $\int dy \, \rho(y) \, K_t(x\,;y) = 1$ for all $t > 0$ and $x \in \Ri$.\hfill$\Box$

\section{Extension properties}

Although $T$ is defined as a group of isometries and $S$ as a contraction semigroup  on $L_\infty$ 
they do not automatically extend  to the $L_p$-spaces.
This requires extra boundedness conditions on the coefficient function $a$ 
and the density function $\rho$.
The following proposition gives necessary and sufficient conditions for  $T$ to extend to a continuous group and
sufficient conditions for $S$ to extend to a continuous semigroup.

\begin{prop}\label{pdvf1.2}
Let $T$ be the group of isometries of $L_\infty(\Ri\,;\rho\,dx)$  defined by
$(T_t\varphi)(y)=\varphi(e^{-tX}y)$.
 The following conditions are equivalent for  all $C\geq 1$ and $\omega\geq0$.
\begin{tabel}
\item\label{pdvf1.2-1}
There is a $p\in[1,\infty\rangle$ such that  $T$ extends to a $($strongly$)$ continuous group on 
$L_p(\Ri\,;\rho \,dx)$ satisfying
the bounds
$\|T_t\|_{p\to p}
\leq C^{1/p}\,e^{\omega |t|/p}$ for all $t\in\Ri$.
\item\label{pdvf1.2-2}
For all $p\in[1,\infty\rangle$ the group $T$ extends to a $($strongly$)$  continuous group on 
$L_p(\Ri\,;\rho \,dx)$ satisfying
the bounds 
$\|T_t\|_{p\to p}
\leq C^{1/p}\,e^{\omega |t|/p}$ for all $t\in\Ri$.
\item\label{pdvf1.2-3}
\hspace{2cm}$\;\;\;\;\;\;a(y)\rho(y)\leq C\,e^{\omega d(x;y)}\,a(x)\rho(x)\;\;\;\;\;\;$ 
\hspace{5mm}for all $x,y\in \Ri$.
\end{tabel}

Moreover, if these conditions are satisfied then the semigroup $S$ extends to a $($strongly$)$ continuous
semigroup on all the $L_p$-spaces, $p\in[1,\infty\rangle$, satisfying
the bounds 
\[
\|S_t\|_{p\to p}
\leq \Big((2\,C)^{1/p}\,e^{\omega^2 t/p}\Big)\wedge \Big(2\,C^{1/p}\,e^{\omega^2 t/p^2}\Big)
\]
if $\omega>0$ and $\|S_t\|_{p\to p}\leq C^{1/p}$ if $\omega=0$, for all $t > 0$.
\end{prop}
\proof\
First assume Condition~\ref{pdvf1.2-1} is satisfied.
Then for all $\varphi\in L_p$ one has
\[
\|T_t\varphi\|^p_p=\int_\Ri dy\,\rho(y)\,|\varphi(e^{-tX}y)|^p
\;\;\;.
\]
Secondly, by a change of variables $x=e^{-tX}y$ one finds
\[
\|T_t\varphi\|^p_p=\int_\Ri dx\,{{a(e^{tX}x)}\over{a(x)}}\,\rho(e^{tX}x)\,|\varphi(x)|^p
      =\int_\Ri dx\,\rho(x)\,\Big({{a(e^{tX}x)\rho(e^{tX}x)}\over{a(x)\rho(x)}}\Big)
\,|\varphi(x)|^p
\;\;\;.
\]
Therefore 
\[
\sup_{x\in\Ri}\Big({{a(e^{tX}x)\rho(e^{tX}x)}\over{a(x)\rho(x)}}\Big)^{1/p}=\|T_t\|_{p\to p}
\leq C^{1/p}e^{\omega |t|/p}
\]
for all $t\in\Ri$ and $x\in \Ri$.
Hence 
\[
a(e^{tX}x)\rho(e^{tX}x)\leq C\,e^{\omega |t|} a(x)\rho(x)
\]
for all $t\in\Ri$ and $x\in \Ri$.
Setting $y=e^{tX}x$ and noting that $d(x\,;y)=|t|$ one deduces that
Condition~\ref{pdvf1.2-3} is satisfied.
Conversely, the same calculation shows that if Condition~\ref{pdvf1.2-3} is satisfied then 
\begin{equation}
\|T_t\varphi\|_p
\leq C^{1/p}e^{\omega |t|/p}\|\varphi\|_p
\label{edvf1.10}
\end{equation}
for all $p\in[1,\infty\rangle$, $\varphi\in L_p$ and
$t\in\Ri$.
In addition if $\varphi\in C_c^\infty$ then one calculates that
\[
\varphi-T_t\varphi=\int^t_0ds\,T_sX_{\rm min}\varphi
\;\;\;.
\]
Hence using (\ref{edvf1.10}) and the density of $C_c^\infty$ in $L_p$
one concludes that $T_t$ extends to a continuous semigroup on $L_p$
satisfying the bounds (\ref{edvf1.10}), i.e., Condition~\ref{pdvf1.2-2} is valid. 
The implication \ref{pdvf1.2-2}$\Rightarrow$\ref{pdvf1.2-3} 
is trivial.

If the conditions are satisfied then $S$ extends to the $L_p$-spaces by (\ref{edvf1.3}).
The estimates on the norms of $S_t$ are established in two steps.
First, if $\omega>0$ then it follows from (\ref{edvf1.3}) and the estimates on $\|T_s\|_{1\to1}$  that 
\[
\|S_t\|_{1\to1}\leq 2\,C\,e^{\omega^2 t}
\]
for all $t > 0$.
Since $S$ is contractive on $L_\infty$ one deduces from interpolation that 
\[
\|S_t\|_{p\to p}
\leq (2\,C)^{1/p}\,e^{\omega^2 t/p}
\]
for all $p\in\langle1,\infty\rangle$ and $t > 0$.
Alternatively, one can reverse the reasoning and use the interpolated bounds 
$\|T_s\|_{p\to p}\leq C^{1/p}\,e^{\omega |s|/p}$ together with  (\ref{edvf1.3})
to calculate that
\[
\|S_t\|_{p\to p}
\leq 2\,C^{1/p}\,e^{\omega^2 t/p^2}
\]
for all $p\in[1,\infty]$ and $t > 0$.

If $\omega=0$ similar arguments apply and both lead to the bounds
$\|S_t\|_{p\to p}\leq C^{1/p}$. \hfill$\Box$

\ruimte

The situation described by the proposition simplifies if $C=1$.
Then Condition~\ref{pdvf1.2-3} together with (\ref{edvf1.61}) implies that
\begin{eqnarray*}
\pm(a\rho)'(y)\,a(y)
&=&\lim_{t\downarrow 0}t^{-1} \Big( (a\rho)(e^{\pm tX}y)-(a\rho)(y) \Big)  \\[5pt]
&\leq&\limsup_{t\downarrow0}t^{-1}(e^{\omega t}-1)(a\rho)(y)=\omega\,(a\rho)(y)
\end{eqnarray*}
for all $y\in \Ri$.
Thus $\|\rho^{-1}(a\rho)'\|_\infty\leq\omega$.
Conversely, if $\|\rho^{-1}(a\rho)'\|_\infty\leq\omega$ then
\[
\rho(e^{tX}y)^{-1}{{d}\over{dt}}\Big(e^{-\omega t}\,(a\rho)(e^{\pm tX}y)\Big)\leq 0
\]
for all $t\geq0$.
Hence Condition~\ref{pdvf1.2-3} is satisfied with $C=1$.
But the condition $\|\rho^{-1}(a\rho)'\|_\infty\leq\omega$
can be expressed in terms of the vector field.
Therefore one has the following corollary.

\begin{cor}\label{cdvf1.3}
 The following conditions are equivalent for  all  $\omega\geq0$.
\begin{tabel}
\item\label{cdvf1.3-1}
There is a $p\in[1,\infty\rangle$ such that  $T$ extends to a continuous group on 
$L_p(\Ri\,;\rho \,dx)$ satisfying
the bounds
$\|T_t\|_{p\to p}
\leq e^{\omega |t|/p}$ for all $t\in\Ri$.
\item\label{cdvf1.3-2}
For all $p\in[1,\infty\rangle$ the group $T$ extends to a continuous group on 
$L_p(\Ri\,;\rho \,dx)$ satisfying
the bounds 
$\|T_t\|_{p\to p}
\leq \,e^{\omega |t|/p}$ for all $t\in\Ri$.
\item\label{cdvf1.3-3}
\hspace{4.5cm}$\|\rho^{-1}(a\rho)'\|_\infty\leq \omega\;\;\;.$ 
\item\label{cdvf1.3-4}
\hspace{3.5cm}$|(\psi,(X+X^*)\varphi)|\leq \omega \, \|\psi\|_q \, \|\varphi\|_p$\hspace{2cm} for 
all $\varphi,\psi\in C^\infty_c(\Ri)$ and for one pair {\rm (}for all pairs{\rm )}
of dual exponents
$p,q\in[1,\infty]$.
\end{tabel}

Moreover, if these conditions are satisfied then the semigroup $S$ extends to a continuous
semigroup on all the $L_p$-spaces, $p\in[1,\infty\rangle$, satisfying
the bounds 
\[
\|S_t\|_{p\to p}
\leq e^{\omega^2 t/p^2}
\]
 for all $t > 0$.
In addition $H_{\rm max}$ satisfies a G{\aa}rding inequality. 
Precisely, 
\[
\RRe (\varphi,H_{\rm max}\varphi)\geq(1-\varepsilon)\|X\varphi\|_2^2-
(4\varepsilon)^{-1})\|X+X^*\|^2_{2\to2}\|\varphi\|_2^2
\]
for all $\varphi\in C_c^\infty(\Ri)$ and $\varepsilon>0$.
\end{cor}
\proof\
The equivalence of the first three conditions and the existence of the extension of the semigroup
$S$ follow from Proposition~\ref{pdvf1.1} and the above discussion.
Conditions~\ref{cdvf1.3-3} and~\ref{cdvf1.3-4} are equivalent because
\begin{eqnarray*}
(\psi,X\varphi)+(X\psi,\varphi)
&=&\int_\Ri dx\,(a\rho)(x) \Big(\psi(x) \, \varphi'(x)
+\psi'(x) \, \varphi(x)\Big)\\[5pt]
&=&\int_\Ri dx\,\rho(x)\Big(\rho(x)^{-1}(a\rho)'(x)\Big)\psi(x) \, \varphi(x)
\end{eqnarray*}
for all $\varphi,\psi\in C_c^\infty(\Ri)$.
It remains to prove the G{\aa}rding inequality.

If $\varepsilon > 0$ then
\begin{eqnarray*}
\RRe (\varphi,H_{\rm max}\varphi)&=&-\RRe(X^*\varphi,X\varphi)\\[5pt]
&=&\|X\varphi\|_2^2-\RRe((X^*+X)\varphi,X\varphi)\\[5pt]
&\geq&\|X\varphi\|_2^2-\|(X^*+X)\varphi\|_2\|X\varphi\|_2\\[5pt]
&\geq&(1-\varepsilon)\|X\varphi\|_2^2-
(4\varepsilon)^{-1} \|X+X^*\|^2_{2\to2}\|\varphi\|_2^2
\end{eqnarray*}
for all $\varphi\in C_c^\infty(\Ri)$.\hfill$\Box$

\ruimte

The corollary, applied with $\omega=0$, gives the following criteria for $T$ or $S$ to extend 
to a  contraction group or semigroup on the $L_p$-spaces.

\begin{prop} \label{pvf37}
The following are equivalent.
\begin{tabel} 
\item \label{pvf37-1}
There is a $p\in[1,\infty\rangle$ such that  $T$ extends to a continuous contraction group on 
$L_p(\Ri\,;\rho \,dx)$.
\item \label{pvf37-2}
For all $p\in[1,\infty\rangle$ the group $T$ extends to a continuous contraction group on 
$L_p(\Ri\,;\rho \,dx)$.
\item \label{pvf37-3}
There is a $p\in[1,\infty\rangle$ such that  $S$ extends to a continuous contraction group on 
$L_p(\Ri\,;\rho \,dx)$.
\item \label{pvf37-4}
For all $p\in[1,\infty\rangle$ the semigroup $S$ extends to a continuous contraction group on 
$L_p(\Ri\,;\rho \,dx)$.
\item \label{pvf37-5}
The function $a\rho$ is constant.
\end{tabel}
\end{prop}
\proof\
The implications \ref{pvf37-5}$\Leftrightarrow$\ref{pvf37-1}$\Leftrightarrow$\ref{pvf37-2}$\Rightarrow$\ref{pvf37-4} 
follow from Corollary~\ref{cdvf1.3} and the implication \ref{pvf37-4}$\Rightarrow$\ref{pvf37-3} is trivial.

The proof of the implication \ref{pvf37-3}$\Rightarrow$\ref{pvf37-5} relies on the reasoning
of Lumer and Phillips.

If Condition~\ref{pvf37-3} is valid for some $p \in [1,2]$ then it follows 
by interpolation with the contraction semigroup on $L_\infty$ that Condition~\ref{pvf37-3} is valid 
for all  $p>2$.
Hence it suffices to show that if  $p \in \langle2,\infty\rangle$ and 
$S$ extends to a continuous contraction group on 
$L_p(\Ri\,;\rho \,dx)$ then the function $a \rho$ is constant, i.e., 
Condition~\ref{pvf37-5} is valid.
Fix $p \in \langle2,\infty\rangle$ and assume $S$ extends to a continuous contraction group on 
$L_p(\Ri\,;\rho \,dx)$.
Then it follows from the Lumer--Phillips theorem, \cite{LuP} Theorem 3.1, that 
the generator $H$ of the semigroup $S$ on $L_p(\Ri\,;\rho \,dx)$ is dissipative. 
So if $[\,\cdot\,,\,\cdot\,]$ is a semi-inner product on $L_p(\Ri\,;\rho \,dx)$ then 
$\RRe [H \varphi, \varphi] \geq 0$ for all $\varphi \in D(H)$.
If $\varphi \in C_c^2(\Ri)$ is real valued then $\varphi \in D(H_{\rm max})$ and 
$H_{\rm max} \varphi \in L_p(\Ri\,;\rho \,dx)$. 
So $\varphi \in D(H)$ and $H_{\rm max} \varphi = H \varphi$.
Moreover, 
\[
\int d( a \, \rho \, \varphi^{p-1}) \, a \, (d \, \varphi)
= \int \rho \, \varphi^{p-1} \,  H_{\rm max} \varphi
= \int \rho \, \varphi^{p-1} \,  H \varphi
= \|\varphi\|_p^{p-2} [H \varphi, \varphi]
\geq 0
\]
where $d = d/dx$.
Hence
\begin{equation}
\int d( a \, \rho \, \varphi^{p-1}) \, a \, (d \, \varphi)
\geq 0
\label{epvf37;1}
\end{equation}
for all real valued $\varphi \in W^{1,\infty}_c(\Ri)$ by approximation.

Next fix $\tau \in C_c^\infty(\Ri)$ such that $0 \leq \tau \leq 1$, $\tau(0) = 1$ and 
$\tau$ is decreasing on $[0,\infty\rangle$.
For all $n \in \Ni$ define $\varphi_n \in W^{1,\infty}_c(\Ri)$ by
\[
\varphi_n 
= (a \rho)^{-1/p} \, (\tau \circ \Phi_n)
\]
where 
\[
\Phi_n(x) 
= n^{-1} \, d(0\,;x)^2
= n^{-1} \Big( \int_0^x a^{-1} \Big)^2
\;\;\; .  \]
Then 
\begin{eqnarray*}
\varphi_n'(x)
& = & - p^{-1} (a \rho)(x)^{-1 - p^{-1}} \, (a \rho)'(x) \, \tau(\Phi_n(x))   \\*[5pt]
& & \hspace{20mm} {}
   + 2 n^{-1} (a \rho)(x)^{-1/p} \, \tau'(\Phi_n(x)) \Big( \int_0^x a^{-1} \Big) a(x)^{-1}
\end{eqnarray*}
and 
\begin{eqnarray*}
(a \rho \, \varphi_n')(x)
& = & - p^{-1} (a \rho)(x)^{-1/p} \, (a \rho)'(x) \, \tau(\Phi_n(x)) \\*[5pt]
& & \hspace{20mm} {}
   + 2 n^{-1} \rho(x) \, (a \rho)(x)^{-1/p} \, \tau'(\Phi_n(x)) \Big( \int_0^x a^{-1} \Big) 
\;\;\; . 
\end{eqnarray*}
Similarly, $(a \rho \, \varphi_n^{p-1})(x) = (a \rho)(x)^{1/p} \, \tau(\Phi_n(x))^{p-1}$
and 
\begin{eqnarray*}
(a \rho \, \varphi_n)'(x)
& = & p^{-1} (a \rho)(x)^{-1 + p^{-1}} \, (a \rho)'(x) \, \tau(\Phi_n(x))^{p-1}  \\*[5pt]
& & \hspace{2mm} {}
   + 2 n^{-1} (p-1) \rho(x) \, (a \rho)(x)^{-1 + p^{-1}} \, \tau(\Phi_n(x))^{p-2} \,
          \tau'(\Phi_n(x)) \Big( \int_0^x a^{-1} \Big) 
\;\;\; . 
\end{eqnarray*}
Then by (\ref{epvf37;1}) it follows that 
\begin{eqnarray*}
0 
& \leq & \int \rho^{-1} d( a \rho \, \varphi_n^{p-1}) \, a \rho \, (d \, \varphi_n)  \\[5pt]
& = & \int dx \bigg( - p^{-2} \rho(x)^{-1} \, (a \rho)(x)^{-1} \, (a \rho)'(x)^2 \Big( \tau(\Phi_n(x)) \Big)^2  \\*[5pt]
& & \hspace{15mm} {} 
    - 2 n^{-1} (1 - 2 p^{-1}) \, (a \rho)(x)^{-1} \, (a \rho)'(x) \, \tau(\Phi_n(x))^{p-1} \,
                            \tau'(\Phi_n(x)) \Big( \int_0^x a^{-1} \Big)    \\*[5pt]
& & \hspace{15mm} {} 
    + 4 n^{-2} (p-1) \rho(x) \, (a \rho)(x)^{-1} \tau(\Phi_n(x))^{p-1} \,
            \Big( \tau'(\Phi_n(x)) \Big)^2 \, d(0\,;x)^2 
      \Bigg)
\;\;\; .
\end{eqnarray*}
Using the estimate $a \, b \leq \varepsilon a^2 + (4 \varepsilon)^{-1} b^2$ for the second term, setting 
$\varepsilon = (2 p (p-2))^{-1}$ and rearranging one finds
\begin{eqnarray}
\lefteqn{
(2 p^2)^{-1} \int \rho^{-1} \, (a \rho)^{-1} \, ((a \rho)')^2 (\tau \circ \Phi_n)^2 
} \hspace{10mm} \nonumber \\*[5pt]
& \leq & n^{-1} \int \rho \, (a \rho)^{-1} \bigg( 
      4 (p-1) (\tau \circ \Phi_n)^{p-2} + 2 (p-2)^2 (\tau \circ \Phi_n)^{2p-2} \bigg) 
      (\tau' \circ \Phi_n) )^2 \, \Phi_n
\hspace*{10mm}
\label{epvf37;2}
\end{eqnarray}
for all $n \in \Ni$.
There are $b,c > 0$ such that 
\[
y \, \Big( 4 (p-1) \tau(y)^{p-2} + 2 (p-2)^2 \tau(y)^{2p-2} \Big) (\tau'(y))^2 \leq c \, e^{-(4b)^{-1} y}
\]
for all $y \in [0,\infty\rangle$.
Then 
\begin{eqnarray*}
\lefteqn{
\bigg( (a \rho)^{-1} \Big( 4 (p-1) (\tau \circ \Phi_n)^{p-2} + 2 (p-2)^2 (\tau \circ \Phi_n)^{2p-2} \Big)
         (\tau' \circ \Phi_n)^2 \, \Phi_n \bigg) (x)
} \hspace{100mm} \\*[5pt]
& \leq & c \, (a \rho)(x)^{-1} \, e^{- d(0;x)^2 (4bn)^{-1}}  \\[5pt]
& = & c \, (4 \pi \, b \, n)^{1/2} \, K_{bn}(0\,;x)
\end{eqnarray*}
uniformly for all $x \in \Ri$ and $n \in \Ni$.
Using Proposition~\ref{pdvf1.21} one deduces that 
\[
\int \rho \, (a \rho)^{-1} \Big( 4 (p-1) (\tau \circ \Phi_n)^{p-2} + 2 (p-2)^2 (\tau \circ \Phi_n)^{2p-2} \Big)
       (\tau' \circ \Phi_n)^2 \, \Phi_n
\leq c \, (4 \pi \, b \, n)^{1/2} 
\]
for all $n \in \Ni$.
Finally  (\ref{epvf37;2}) and the monotone convergence theorem  establishes that 
\begin{eqnarray*}
(2 p^2)^{-1} \int \rho^{-1} \, (a \rho)^{-1} \Big( (a \rho)' \Big)^2 
& = & \lim_{n \to \infty} (2 p^2)^{-1} \int \rho^{-1} \, (a \rho)^{-1} \, \Big( (a \rho)' \Big)^2 
         (\tau \circ \Phi_n)^2  \\[5pt]
& \leq & \lim_{n \to \infty} n^{-1} \, c \, (4 \pi \, b \, n)^{1/2} 
= 0
\;\;\; .
\end{eqnarray*}
Therefore $(a \rho)' = 0$ as required.\hfill$\Box$

\ruimte

In the unweighted case, i.e., $\rho=1$, the proposition establishes that $S$ extends to  a contraction semigroups on  one of the $L_p$-spaces with $p<\infty$  only in the case that $X$ is proportional
to $d/dx$.

\section{Examples}\label{Svf4}

Next we give two examples of rather unexpected properties
although  there is nothing inherently pathological about the weight $\rho$ or the coefficient $a$.
In fact in both examples $\rho=1$ and the coefficient $a$ of the vector field is strictly positive,  smooth and uniformly bounded.
The first example gives a continuous group $T$ and semigroup $S$ which do not extend from $L_\infty$ to the other  $L_p$ spaces.
The principal reason for this singular behaviour is the fact that $\inf a = 0$, i.e., there is a mild degeneracy at infinity.

\begin{exam} \label{xtvf1}
Let $\rho = 1$.
For all $n \in \Ni_0$ define $h_n = n!^{-1}$.
Define $y_n \in \Ri$ for all $n \in \Ni_0$ by $y_0 = 0$ and inductively
\[
y_{n+1} = y_n + 4^{-1} (h_n + h_{n+1}) + 2^{-1}
\]
for all $n \in \Ni$.
Define $\tilde a \colon \Ri \to \langle0,\infty\rangle$ by
\[
\tilde a(x)
= \left\{ \begin{array}{ll}
   h_n & \mbox{if } x \in [y_n - 4^{-1} h_n, y_n + 4^{-1} h_n \rangle \;\;\;\; (n \in \Ni_0) \;\;\; , \\[5pt]
   1   & \mbox{if } x \in [y_n + 4^{-1} h_n, y_n + 4^{-1} h_n + 2^{-1} \rangle \;\;\;\; (n \in \Ni_0)  \;\;\; , \\[5pt]
   1   & \mbox{if } x \in \langle-\infty,0]  \;\;\; . 
           \end{array} \right.
 \]
Then $\tilde a(y_{n}) = h_n$ and $\int_{y_n}^{y_{n+1}} dx \, \tilde a(x)^{-1} = 1$
for all $n \in \Ni$.
Next we regularize $\tilde a^{-1}$.
For all $n \in \Ni_0$ let $\chi_n \in C_c^\infty(\Ri)$ be such that 
$\chi_n \geq 0$, $\int \chi_n = 1$, $\supp \chi_n \subseteq [-8 ^{-1}h_n, 8^{-1} h_n]$
and $\chi_n(-x) = \chi_n(x)$ for all $x \in \Ri$.
Define $a \in C^\infty(\Ri)$ by
\[
a(x)^{-1} 
= \left\{ \begin{array}{ll}
   (\chi_0 * \tilde a^{-1})(x) & \mbox{if } x \leq 0 \;\;\; , \\[5pt]
   (\chi_n * \tilde a^{-1})(x) & \mbox{if } n \in \Ni_0 \mbox{ and }
          x \in [y_{n} - 4^{-1} h_{n} - 4^{-1}, y_{n} + 4^{-1} h_{n} + 4^{-1} \rangle \;\;\; .
           \end{array} \right.
\]
Then $a(y) = h_n$ for all $y\in [y_n-8 ^{-1}h_n, y_n+8 ^{-1}h_n]$
and $\int_{y_n}^{y_{n+1}} dx \, a(x)^{-1} = 1$
for all $n \in \Ni$.
Hence $d(y_n\,;y_{n+1}) = 1$ for all $n \in \Ni$.
But $a(y_{n}) = (n+1) \, a(y_{n+1})$ for all $n \in \Ni$.
Therefore Condition~\ref{pdvf1.2-3} of Proposition~\ref{pdvf1.2} is not valid.
In particular the group $T$ does not extend to any of the other $L_p$ spaces.
Next we show that the semigroup $S$ also does not extend to another $L_p$ space.

Let $p \in [1,\infty\rangle$, $t > 0$ and let $q$ be the dual exponent of $p$.
For all $n \in \Ni$ set $I_n = [y_n - 8^{-1} h_n, y_n + 8^{-1} h_n]$. 
Let $n \in \Ni$.
Set $\varphi = \one_{I_{n+1}}$ and $\psi = \one_{I_n}$.
Then $\|\varphi\|_p = |I_{n+1}|^{1/p}$ and $\|\psi\|_q = |I_n|^{1/q}$.
Moreover, 
\begin{eqnarray*}
(\psi, S_t \varphi)
& = & (4\pi t)^{-1/2} \int_{I_n} dx \int_{I_{n+1}} dy \, a(y)^{-1} \, e^{- d(x;y)^2 (4t)^{-1}}  \\[5pt]
& \geq & (4\pi t)^{-1/2} \int_{I_n} dx \int_{I_{n+1}} dy \, a(y)^{-1} \, e^{- 3 d(x;y)^2 t^{-1}}  \\[5pt]
& = & (4\pi t)^{-1/2} |I_n| \, |I_{n+1}| \, h_{n+1}^{-1}  \, e^{- 3 d(x;y)^2 t^{-1}}
\;\;\; .  
\end{eqnarray*}
So 
\[
\|S_t\|_{p \to p}
\geq (4\pi t)^{-1/2} |I_n|^{1/p} \, |I_{n+1}|^{1/q} \, h_{n+1}^{-1} \, e^{- 3 d(x;y)^2 t^{-1}} 
= (64\pi t)^{-1/2} (n+1)^{1/p}
\;\;\; .  \]
Hence the operator $S_t$ on $L_\infty$ does not extend to a continuous operator 
on $L_p$ for any $p\in[1,\infty\rangle$ or $t > 0$.\hfill$\Box$
\end{exam}

In the next example the coefficient $a$ of $X$ is uniformly bounded above and below by a positive constant
but $\sup a'=\infty$
The semigroup
$S$ extends to a continuous semigroup on all the $L_p$-spaces but
 the real part of the generator of $S$ on $L_2$ is not lower semibounded.
This contrasts with the case of continuous self-adjoint semigroups where boundedness of the semigroup
immediately implies lower semiboundedness of the generator.

\begin{exam}\label{xdvf1}
First, let $\rho=1$ and let $\chi\in C_c^\infty(\Ri)$ be such that $0\leq \chi\leq 3$, $\chi'\geq0$, $\chi(x)=0$
if  $x\leq 0$, $\chi(x)=3$ if $x \geq 3$ and $\chi(x)=x$ if $1\leq x\leq2$.
Define $a\colon \Ri\to[1,4]$ by
\[
a(x)=1+\sum_{n=1}^\infty\Big(\chi(n(x-16n))-\chi(n(x-(16n+8))\Big)
\;\;\;.
\]
Thus $a=1$ on an infinite sequence of intervals of length almost equal to $8$ spaced at distance $8$ one from the other.
On the intermediate intervals $a$ increases smoothly to the value $4$ and then decreases in a similar fashion to the value $1$.
The rate of increase and decrease, however,  becomes larger with the distance of the interval from the origin.
Nevertheless $a\in C^\infty(\Ri)$ and the bounds of Proposition~\ref{pdvf1.2}.\ref{pdvf1.2-3} are valid 
with $C = 4$ and $\omega = 0$.
In particular $S_t$ extends to the $L_p$-spaces and $\|S_t\|_{p\to p}\leq 4^{1/p}$.

Secondly, let $n\in\Ni$ with $n \geq 4$.
Let $\psi \in C^\infty(\Ri)$ be such that $\psi(x) = 3$ for all $x\leq 16n+8$,
$0\leq \psi'\leq n^{1/2}$, $\psi'(x)=0$ for all $x\geq 16n+8+4 n^{-1}$ and 
$\psi'(x)=n^{1/2}$ for all $x\in[16n+8+ n^{-1},16n+8+2 n^{-1}]$.
Then $3 \leq \psi(16n+8+4 n^{-1}) \leq 5$.
Now define $\varphi \in C_c^\infty(\Ri)$ by
\[
\varphi(x)
= \left\{
\begin{array}{ll}
\chi(x-(16n+4)) & \mbox{if }x\leq 16n+8\\[5pt]
\psi(x)  & \hspace*{-15mm} \mbox{if }x\in[16n+8,16n+8+4 n^{-1}]\\[5pt]
3^{-1} \psi(16n+8+4 n^{-1}) \Big( 3 -\chi(x-(16n+8+4 n^{-1})\Big) & \mbox{if }x\geq 16n+8+4 n^{-1}
\end{array}
\right.
\]
Then $\|\varphi\|_2 \leq 5 \cdot (12)^{1/2} = (300)^{1/2}$ and 
\[
\|\varphi'\|_2
\leq 2 \|\chi'\|_\infty + n^{1/2}(4 n^{-1})^{1/2} + 3^{-1}\psi(16n+8+4 n^{-1})\|\chi'\|_\infty
\leq 2 + 4 \|\chi'\|_\infty
\;\;\;.
\]
But $a' \, a \, \varphi \, \varphi'\leq 0$ and 
\[
-(a'\varphi,X\varphi)
\geq \int^{16n+8+2 n^{-1}}_{16n+8+ n^{-1}} (-a' \, a \, \varphi \, \varphi')
\geq \int^{16n+8+2 n^{-1}}_{16n+8+ n^{-1}} n \cdot 2 \cdot 3 \cdot n^{1/2}
= 6 n^{1/2}
\]
by the previous estimates.
Therefore 
\begin{eqnarray*}
\RRe(\varphi,H_{\rm min} \varphi)
& = & \|X\varphi\|^2_2 + \RRe(a'\varphi,X\varphi)\\[5pt]
& \leq & \|a\|_\infty^2 (2+4\|\chi'\|_\infty)^2 - 8n^{1/2}
\leq -300^{-1} \Big(6 n^{1/2} - 16(2+ 4 \|\chi'\|_\infty)^2 \Big) \|\varphi\|_2^2
\;\;\;.
\end{eqnarray*}
Consequently, $\RRe H_{\rm min}$ is not lower semibounded.
This is despite the uniform boundedness of $S$ on $L_2$.

Next, since $S$  is uniformly bounded on each of the $L_p$-spaces, the spectrum
 $\sigma(H)$ of the generator $H$ of the semigroup on $L_p$ is contained in the right half-plane.
But  $a(x) \in [1,4]$ for all $x \in \Ri$. 
Therefore
$4^{-1} |x-y| \leq d(x\,;y) \leq |x-y|$
and Proposition~\ref{pdvf1.21} implies that 
\[
K_t(x\,;y) 
\leq  (4\pi t)^{-1/2} \, e^{-|x-y|^2 (64 t)^{-1}}
\]
for all $x,y \in \Ri$ and $t > 0$.
Hence it follows from \cite{Kun} or \cite{LisV} that  $\sigma(H)$
is independent of $p \in [1,\infty]$.
On the other hand $\RRe H_{\rm min}$ is not lower semibounded on $L_2$ and the 
above estimates establish  that $\langle-\infty,0] \subset \Theta(H)$, the $L_2$-numerical range 
of $H$.
Therefore $\Theta(H)\neq \sigma(H)$ on $L_2$.

In fact this example illustrates the extreme situation that the spectrum of $H$ is contained in the right half plane but the numerical range is the whole complex plane.
This follows since one can establish that the  numerical range $\Theta(H) = \Ci$
by a small modification  of the foregoing estimates applied to the function
$\tilde \varphi \in C_c^\infty(\Ri)$ defined by
\[
\tilde \varphi(x) = e^{i \lambda x} \, \tau(x) + \varphi(x)
\;\;\;,
\]
where $\lambda \in \Ri$ and $\tau \in C_c^\infty(\langle-1,4\rangle)$ is fixed such that 
$0 \leq \tau \leq 1$ and $\tau|_{[0,3]} = 1$.
One also uses  the observation that  the numerical range is convex.

Finally note that the semigroup $S$ has a bounded holomorphic extension to the open right half-plane
on each of the $L_p$-spaces, $p\in[1,\infty \rangle$.
This follows from the explicit form of the kernel given in Propositions~\ref{pdvf1.21}.
Therefore the operator $H$ is of type $S_{0+}$.
Nevertheless, since $\Theta(H) = \Ci$ the operator $H$ is not sectorial.\hfill$\Box$
\end{exam}

\section{Volume doubling}\label{Svf5}

Let $V(x\,;r)$ denote the measure of the ball of radius $r$ centred at $x$,
i.e., the set $\{y:d(x\,;y)<r\}=\langle e^{-rX}x,e^{rX}x\rangle$.
Then $V$ is defined, as usual, to have the volume doubling property if 
there is a $c>0$ such that 
\[
V(x\,;2r)\leq c\,V(x\,;r)
\]
for all $r>0$.
This property can be immediately related to the conditions of Proposition~\ref{pdvf1.2}
which are necessary and sufficient for the continuous extension of $T$ to the $L_p$-spaces.

\begin{prop}\label{pvf5.1}
\mbox{}
\begin{tabel}
\item\label{pvf5.1-1}
If the equivalent conditions of  Proposition~$\ref{pdvf1.2}$ are satisfied
then
\begin{equation}
V(x\,;2r)\leq 2\,C^2\,e^{3\omega}\,V(x\,;r)
\label{evf5.1}
\end{equation}
for all $x\in\Ri$ and $r\in\langle0,1]$
where $C$ and $\omega$ are the parameters of Proposition~$\ref{pdvf1.2}$.
Moreover if $\omega=0$ then $(\ref{evf5.1})$ is valid for all $x\in\Ri$ and $r>0$.
\item\label{pvf5.1-2}
If there exist $c > 0$ and a function $v \colon \langle 0,\infty \rangle \to \Ri$ such that 
\[
c^{-1} \, v(r) \leq V(x\,;r) \leq c \, v(r)
\]
for all $x \in \Ri$ and $r \in \langle0,1]$
then Condition~{\rm \ref{pdvf1.2-3}} of Proposition~$\ref{pdvf1.2}$ is satisfied with $\omega = 0$.
\end{tabel}
\end{prop}
\proof\
It follows by definition that 
\[
V(x\,;r) =\int_{e^{-rX}x}^{e^{rX}x}dy\,\rho(y)
\;\;\;.
\]
But 
\[
{{d}\over{dr}}V(x\,;r)=(a\rho)(e^{rX}x)+(a\rho)(e^{-rX}x)
\;\;\;.
\]
Hence
\[
V(x\,;r)=\int^r_0 ds\,\Big((a\rho)(e^{sX}x)+(a\rho)(e^{-sX}x)\Big)=
\int^r_{-r} ds\,(a\rho)(e^{sX}x)
\;\;\;.
\]
Therefore if Condition~\ref{pdvf1.2-3} of Proposition~\ref{pdvf1.2}
is satisfied one estimates that
\[
2\,C^{-1}r\,e^{-\omega r}(a\rho)(x)\leq V(x\,;r)\leq 2\,C\,r\,e^{\omega r}(a\rho)(x)
\]
for all $x\in \Ri$ and $r>0$.
These bounds imply (\ref{evf5.1}) for all $x\in\Ri$ and $r\in\langle0,1]$
or, if $\omega=0$, for all $r>0$.

If, however, the assumptions of the second statement are valid
then 
\[
c^{-1} \, v(r) 
\leq V(x\,;r) 
= \int_0^r ds \, (a\rho)(e^{sX}x) + (a\rho)(e^{-sX}x)
\leq r \max_{y \in [e^{-X} x, e^X x]} (a\rho)(y)
\]
for all $x \in \Ri$ and $r \in \langle0,1]$.
Similarly
\[
c \, v(r) \geq r \min_{y \in [e^{-X} x, e^X x]} (a\rho)(y)
\;\;\;.  \]
Hence  there exists a $c_1 > 0$ such that 
$c_1^{-1} \, r \leq v(r) \leq c_1 \, r$ for all $r \in \langle0,1]$.
But then 
\begin{eqnarray*}
2 (a\rho)(x) 
& = & \lim_{r \downarrow 0} r^{-1} \int_0^r ds \, (a\rho)(e^{sX}x) + (a\rho)(e^{-sX}x)  \\[5pt]
& = & \lim_{r \downarrow 0} r^{-1} \, V(x\,;r)
\leq \limsup_{r \downarrow 0} r^{-1} \, c \, v(r)
\leq c \, c_1
\end{eqnarray*}
for all $x \in \Ri$.
Similarly $2 (a\rho)(x) \geq (c \, c_1)^{-1}$.
Hence $(2 c \, c_1)^{-1} \leq a \rho \leq 2^{-1} c \, c_1$ and 
Condition~\ref{pdvf1.2-3} of Proposition~\ref{pdvf1.2} is satisfied with $\omega = 0$.
\hfill$\Box$

\subsection*{Acknowledgement}
This work was completed whilst the second named author was a guest of the Department of Mathematics at the University of Auckland.

\end{document}